\documentclass[preprints,article,accept,pdftex,moreauthors]{Definitions/mdpi} 
\firstpage{1} 
\pubvolume{1}
\issuenum{1}
\articlenumber{0}
\pubyear{2023}
\copyrightyear{2023}
\datereceived{ } 
\daterevised{ } 
\dateaccepted{ } 
\datepublished{ } 
\hreflink{https://doi.org/} 

\newcommand{\Unif}{\mathrm{Unif}}

\newcommand{\Var}{\mathrm{Var}}
\newcommand{\Cov}{\mathrm{Cov}}

\catcode`~=11 \def\UrlSpecials{\do\~{\kern -.15em\lower .7ex\hbox{~}\kern .04em}} \catcode`~=13 

\allowdisplaybreaks[1]

\newcommand{\sym}{\mathrm{s}}

\newcommand{\calC}{\mathcal{C}}

\newcommand{\calU}{\mathcal{U}}

\newcommand{\calX}{\mathcal{X}}

\theoremstyle{plain}
\newtheorem{thm}{\protect\theoremname}
\theoremstyle{plain}
\newtheorem{conjecture}{\protect\conjecturename}
\theoremstyle{plain}
\newtheorem{prop}{\protect\propositionname}
\theoremstyle{plain}
\newtheorem{cor}{\protect\corollaryname}
\theoremstyle{plain}
\newtheorem{lem}{\protect\lemmaname}

\providecommand{\corollaryname}{Corollary}

\providecommand{\lemmaname}{Lemma}
\providecommand{\propositionname}{Proposition}

\providecommand{\theoremname}{Theorem}
\providecommand{\conjecturename}{Conjecture}

\Title{Dimension-Free Bounds for the Union-Closed Sets Conjecture}

\TitleCitation{Dimension-Free Bounds for the Union-Closed Sets Conjecture}

\Author{Lei Yu $^{1}$\orcidA{}}

\AuthorNames{Lei Yu}

\AuthorCitation{Yu, L.}

\address{%
$^{1}$ \quad School of Statistics and Data Science, LPMC, KLMDASR,
and LEBPS, Nankai University,  Tianjin 300071, China; leiyu@nankai.edu.cn}

\abstract{The union-closed sets conjecture states that in any nonempty union-closed
family $\mathcal{F}$ of subsets of a finite set, there exists an
element contained in at least a proportion $1/2$ of the sets of $\mathcal{F}$.
Using the information-theoretic method, Gilmer \cite{gilmer2022constant}
recently showed that there exists an element contained in at least
a proportion $0.01$ of the sets of such $\mathcal{F}$. He conjectured
that his technique can be pushed to the constant $\frac{3-\sqrt{5}}{2}$
which was subsequently confirmed by several researchers \cite{sawin2022improved,chase2022approximate,alweiss2022improved,pebody2022extension}.
Furthermore, Sawin \cite{sawin2022improved} showed that Gilmer's
technique can be improved to obtain a bound better than $\frac{3-\sqrt{5}}{2}$, but this new bound is not explicitly given by Sawin. 
This paper further improves Gilmer's technique to derive new bounds
in the optimization form for the union-closed sets conjecture. These
bounds include Sawin's improvement as a special case. By providing
cardinality bounds on auxiliary random variables, we make Sawin's
improvement computable, and then evaluate it numerically which yields
a bound around $0.38234$, slightly better than $\frac{3-\sqrt{5}}{2}\approx0.38197$. }


\begin{document}

\section{Introduction}

This paper concerns the union-closed conjecture which is described
in the information-theoretic language as follows. Note that each set
$B\subseteq[n]:=\{1,2,\dots,n\}$ uniquely corresponds to an $n$-length
sequence $x^{n}:=(x_{1},x_{2},\dots,x_{n})\in\Omega^{n}$ with $\Omega:=\{0,1\}$
in the way that $x_{i}=1$ if $i\in B$ and $x_{i}=0$ otherwise.
So, a family $\mathcal{F}$ of subsets of $[n]$ uniquely corresponds
to a subset $A\subseteq\Omega^{n}$. Denote the (element-wise) OR
operation for two finite $\Omega$-valued sequences as $x^{n}\vee y^{n}:=(x_{i}\vee y_{i})_{i\in[n]}$
with $x^{n},y^{n}\in\Omega^{n}$, where $\vee$ is the OR operation.
The family $\mathcal{F}$ is closed under the union operation (i.e.,
$F\cup G\in\mathcal{F},\forall F,G\in\mathcal{F}$) if and only if
the corresponding set $A\subseteq\Omega^{n}$ is closed under the
OR operation (i.e., $x^{n}\vee y^{n}\in A,\forall x^{n},y^{n}\in A$).

Let $A\subseteq\Omega^{n}$ be closed under the OR operation. Let
$X^{n}:=(X_{1},X_{2},\dots,X_{n})$ be a random vector uniformly distributed
on $A$, and denote $P_{X^{n}}=\Unif(A)$ as its distribution (or
probability mass function, PMF). We are interested in estimating 
\[
p_{A}:=\max_{i\in[n]}P_{X_{i}}(1)
\]
where $P_{X_{i}}$ is the distribution of $X_{i}$, and hence, $P_{X_{i}}(1)$
is the proportion of the sets containing the element $i$ among all
sets in $\mathcal{F}$. Frankl made the following conjecture.
\begin{conjecture}[Frankl Union-Closed Sets Conjecture]
$p_{A}\ge1/2$ for any OR-closed set $A$.
\end{conjecture}
This conjecture equivalently states that for any union-closed family
$\mathcal{F}$, there exists an element contained in at least a proportion
$1/2$ of the sets of $\mathcal{F}$. Since the union-closed conjecture
was posed by Peter Frankl in 1979, it had attracted a great deal of
research interest; see, e.g., \cite{balla2013union,johnson1998union,karpas2017two,knill1994graph,wojcik1999union}.
We refer readers to the survey paper \cite{bruhn2015journey} for
more details. Gilmer \cite{gilmer2022constant} made a breakthrough
recently, showing that this conjecture holds with constant $0.01$.
Gilmer's method used a clever idea from information theory in which
two independent random vectors were constructed. It was conjectured
by him that his method can improve the constant to $\frac{3-\sqrt{5}}{2}$,
which is now confirmed by several groups of researchers \cite{alweiss2022improved,chase2022approximate,sawin2022improved,pebody2022extension}.
This constant is shown to be the best for an approximate version of
the union-closed sets problem \cite{chase2022approximate}. Moreover,
Sawin \cite{sawin2022improved} further develops Gilmer's idea by
allowing the two random vectors to depend with each other. Such a
technique was in fact used by the present author in several existing
works \cite{YuTan2020_exact,yu2021strong_article,yu2022exact}. By
this technique, Sawin \cite{sawin2022improved} showed that the constant
can be improved to a value that is strictly larger than $\frac{3-\sqrt{5}}{2}$.
However, without cardinality bounds on auxiliary random variables,
Sawin's constant is difficult to compute, and hence, the accurate
value of this improved constant is not explicitly given in \cite{sawin2022improved}.

The present paper further develops Gilmer's (or Sawin's) technique to derive
new constants (or bounds) in the optimization form for the union-closed
sets conjecture. These bounds include Sawin's improvement as a special
case. By providing cardinality bounds on auxiliary random variables,
we make Sawin's improvement computable, and then evaluate it numerically
which yields a bound around $0.38234$, slightly better than $\frac{3-\sqrt{5}}{2}\approx0.38197$.


\section{Main Results}

To state our result, we need to introduce some notations. Since we only
consider distributions on finite alphabets, we do not distinguish
between the terms ``distributions'' and ``probability mass functions''.
For a pair of distributions $(P_{X},P_{Y})$, a coupling of $(P_{X},P_{Y})$
is a joint distribution $P_{XY}$ whose marginals are respectively
$P_{X},P_{Y}$. For a distribution $P_{X}$ defined on a finite alphabet
$\calX$, a coupling $P_{XX'}$ of $(P_{X},P_{X})$ is called symmetric
if $P_{XX'}(x,y)=P_{XX'}(y,x)$ for all $x,y\in\calX$. Denote $\calC_{\sym}(P_{X})$
as the set of symmetric couplings of $(P_{X},P_{X})$. Denote $\delta_{x}$
as the Dirac measure with atom at $x$.

For a joint distribution $P_{XY}$, the (Pearson) correlation coefficient
between $(X,Y)\sim P_{XY}$ is defined by 
\[
\rho_{\mathrm{p}}(X;Y):=\left\{ \begin{array}{ll}
\frac{\Cov(X,Y)}{\sqrt{\Var(X)\Var(Y)}}, & \Var(X)\Var(Y)>0\\
0, & \Var(X)\Var(Y)=0
\end{array}\right..
\]
The maximal correlation between $(X,Y)\sim P_{XY}$ is defined by
\begin{align*}
\rho_{\mathrm{m}}(X;Y) & :=\rho_{\mathrm{p}}(f(X);g(Y))\\
 & =\sup_{f,g}\left\{ \begin{array}{ll}
\frac{\Cov(f(X),g(Y))}{\sqrt{\Var(f(X))\Var(g(Y))}}, & \Var(f(X))\Var(g(Y))>0\\
0, & \Var(f(X))\Var(g(Y))=0
\end{array}\right.,
\end{align*}
where the supremum is taken over all  pairs of real-valued functions
$\left(f,g\right)$  such that $\Var(f(X))\Var(g(Y))<\infty$. Note
that $\rho_{\mathrm{m}}(X;Y)\in[0,1]$, and moreover, $\rho_{\mathrm{m}}(X;Y)=0$
if and only if $X,Y$ are independent. Moreover, $\rho_{\mathrm{m}}(X;Y)$
is equal to the second largest singular value of the matrix $\left[\frac{P_{XY}(x,y)}{\sqrt{P_{X}(x)P_{Y}(y)}}\right]_{(x,y)}$;
see, e.g., \cite{witsenhausen1975sequences}. Clearly, the largest
singular value of the matrix $\left[\frac{P_{XY}(x,y)}{\sqrt{P_{X}(x)P_{Y}(y)}}\right]_{(x,y)}$
is equal to  $1$ with corresponding eigenvectors $(\sqrt{P_{X}(x)})_{x}$ and
$(\sqrt{P_{Y}(y)})_{y}$.

Denote for $p,q,\rho\in[0,1]$, 
\begin{align*}
z_{1} & :=pq-\rho\sqrt{p(1-p)q(1-q)}\\
z_{2} & :=pq+\rho\sqrt{p(1-p)q(1-q)}
\end{align*}
and 
\begin{equation}
\varphi(\rho,p,q):=\mathrm{median}\left\{ \max\{p,q,p+q-z_{2}\},1/2,\min\{p+q,p+q-z_{1}\}\right\} ,\label{eq:-4}
\end{equation}
where $\mathrm{median}A$ denotes the median value of elements in a multiset $A$. We regard  the set in \eqref{eq:-4} as a multiset which means $\mathrm{median}\{a,a,b\}=a$.  
Denote $h(a)=-a\log_{2}a-(1-a)\log_{2}(1-a)$ for $a\in[0,1]$ as
the binary entropy function. Define for $t>0$, 
\begin{equation}
\Gamma(t):=\sup_{P_{\rho}}\inf_{P_{p}:\mathbb{E}h(p)>0,\mathbb{E}p\le t}\mathbb{E}_{\rho}\left[\inf_{P_{pq}\in\calC_{\sym}(P_{p}):\rho_{\mathrm{m}}(p;q)\le\rho}\frac{\mathbb{E}_{p,q}h(\varphi(\rho,p,q))}{\mathbb{E}h(p)}\right],\label{eq:-5}
\end{equation}
where the supremum over $P_{\rho}$ and the infimum over $P_{p}$
are both taken over all finitely supported probability distributions
on $[0,1]$. 

Our main results are as follows.
\begin{thm}
\label{thm:If-given-,} If $\Gamma(t)>1$ for some $t\in(0,1/2)$,
then $p_{A}\ge t$ for any OR-closed $A\subseteq\Omega^{n}$ (i.e.,
for any union-closed family $\mathcal{F}$, there exists an element
contained in at least a proportion $t$ of the sets of $\mathcal{F}$).
\end{thm}
The proof of Theorem \ref{thm:If-given-,} is given in Section 2 by
using a technique based on coupling and entropy. It is essentially
the same as the technique used by Sawin \cite{sawin2022improved}. However,
prior to Sawin's work, such a technique was used by the present author
in several works; see \cite{YuTan2020_exact,yu2021strong_article,yu2022exact}.

Equivalently, Theorem \ref{thm:If-given-,} states that $p_{A}\ge t_{\max}$
for any OR-closed $A\subseteq\Omega^{n}$, where $t_{\max}:=\sup\{t\in(0,1/2):\Gamma(t)>1\}.$
To compute $\Gamma(t)$ or its lower bounds numerically, it requires
to upper bound the cardinality of the support of $P_{p}$ in the outer
infimum in \eqref{eq:-5}, since otherwise, infinitely many parameters
are needed to optimize. This is left to be done in a future work.
The following gives a computable bound.

{If we choose $P_{\rho}=\delta_{0}$, then Theorem
\ref{thm:If-given-,} implies Gilmer's bound in \cite{gilmer2022constant},
since for this case, the couplings constructed in the proof of Theorem
\ref{thm:If-given-,} (given in the next section) turn to be independent,
coinciding with Gilmer's construction. On the other hand, if we choose
$P_{\rho}=\delta_{1}$, then the couplings constructed in our proof
are arbitrary.  In fact, we can make a choice of $P_{\rho}$ better
than these two special cases. As suggested by Sawin \cite{sawin2022improved},
we can choose $P_{\rho}=(1-\alpha)\delta_{0}+\alpha\delta_{1}$ which
in fact leads to an optimization over mixtures of independent couplings
and arbitrary couplings. This final choice yields the following bound. }

{Substituting $\rho=0$ and $1$ respectively into
$\varphi(\rho,p,q)$ yields 
\begin{align}
\varphi(0,p,q) & =p+q-pq,\\
\varphi(1,p,q) & =\mathrm{median}\left\{ \max\{p,q\},1/2,p+q\right\} ,\label{eq:-4-1}
\end{align}
where in the evaluation of $\varphi(1,p,q),$ the following facts
were used: 1) 
\[
p+q-pq-\sqrt{p(1-p)q(1-q)}\le\max\{p,q\}
\]
for all $p,q\in[0,1]$; 2) if $p+q\le1$, then 
\[
p+q-pq+\sqrt{p(1-p)q(1-q)}\ge p+q,
\]
and otherwise, 
\[
1/2<\max\{p,q\}\le p+q-pq+\sqrt{p(1-p)q(1-q)}.
\]
By defining 
\begin{align*}
g(P_{pq},\alpha) & :=(1-\alpha)\mathbb{E}_{(p,q)\sim P_{p}^{\otimes2}}h(p+q-pq)\\
 & \qquad+\alpha\mathbb{E}_{(p,q)\sim P_{pq}}h(\varphi(1,p,q))
\end{align*}
and substituting $P_{\rho}=(1-\alpha)\delta_{0}+\alpha\delta_{1}$
into Theorem \ref{thm:If-given-,}, one obtains the following simpler
bound. }
\begin{prop}
\label{prop:For-,where-the}For $t\in(0,1/2)$, 
\begin{equation}
\Gamma(t)\ge\hat{\Gamma}(t):=\sup_{\alpha\in[0,1]}\inf_{\textrm{symmetric }P_{pq}:\mathbb{E}h(p)>0}\frac{g(P_{pq},\alpha)}{\mathbb{E}h(p)},\label{eq:-9}
\end{equation}
where the infimum is taken over all distributions $P_{pq}$ of the form
$(1-\beta)Q_{a_{1},a_{2}}+\beta Q_{b_{1},b_{2}}$ with 
\begin{equation}
0\le a:=\frac{a_{1}+a_{2}}{2}\le t<b:=\frac{b_{1}+b_{2}}{2}\le1\label{eq:ab}
\end{equation}
and $\beta=0$ or $\beta=\frac{t-a}{b-a}>0$ such that\footnote{Note that $\mathbb{E}h(p)=0$ if and only if $P_{pq}$ is a convex
combination of $\delta_{(0,0)}$, $\delta_{(0,1)}$, $\delta_{(1,0)}$,
and $\delta_{(1,1)}$.} $\mathbb{E}h(p)>0$. Here, 
\begin{equation}
Q_{x,y}:=\frac{1}{2}\delta_{(x,y)}+\frac{1}{2}\delta_{(y,x)}\label{eq:Q}
\end{equation}
 with $\delta_{(x,y)}$ denoting the Dirac measure at $(x,y)$.
\end{prop}
As a consequence of two results above, we have the following corollary.
\begin{cor}
\label{cor:If-given-,-1} If $\hat{\Gamma}(t)>1$ for some $t\in(0,1/2)$,
then $p_{A}\ge t$ for any OR-closed $A\subseteq\Omega^{n}$.
\end{cor}
The proof of Corollary \ref{thm:If-given-,} is given in Section 3.

The lower bound in \eqref{eq:-9} without the cardinality bound on
the support of $P_{pq}$ was given by Sawin \cite{sawin2022improved},
which was used to show $p_{A}>\frac{3-\sqrt{5}}{2}.$ However, thanks
to the cardinality bound, we can numerically compute the best bound
on $p_{A}$ that can be derived using $\hat{\Gamma}(t)$. That is,
$p_{A}\ge\hat{t}_{\max}$ for any OR-closed $A\subseteq\Omega^{n}$,
where $\hat{t}_{\max}:=\sup\{t\in(0,1/2):\hat{\Gamma}(t)>1\}.$ Numerical
results\footnote{\label{fn:Our-code}Our code can be found  on the author's homepage https://leiyudotscholar.wordpress.com/}
show that if we set $\alpha=0.035,t=0.38234$, then the optimal $P_{pq}=(1-\beta)Q_{a,a}+\beta Q_{a,1}$
with $a\approx0.3300622$ and $\beta\approx0.1560676$ which leads
to the lower bound $\hat{\Gamma}(t)\ge1.00000889$. Hence, $p_{A}\ge0.38234$
for any OR-closed $A\subseteq\Omega^{n}$. This is slightly better
than the previous bound $\frac{3-\sqrt{5}}{2}\approx0.38197$. The
choice of $(\alpha,t)$ in our evaluation is nearly optimal.  {More
decimal places of Sawin's bound (or equivalently, $\hat{t}_{\max}$)
were computed by Cambie in \cite{cambie2022better}, i.e., $0.382345533366702\le\hat{t}_{\max}\le0.382345533366703$
which is attained by the choice $\alpha\approx0.03560698136437784$.
 This more precise evaluation can be also verified using our code
in Footnote \ref{fn:Our-code}. }

\section{Proof of Theorem \ref{thm:If-given-,}}

Denote $H(X)=-\sum_{x}P_{X}(x)\log P_{X}(x)$ as the Shannon entropy
of a random variable $X\sim P_{X}$. Let $A\subseteq\Omega^{n}$ be
closed under the OR operation. We assume $|A|\ge2$. This is because,
Theorem \ref{thm:If-given-,} holds obviously for singletons $A$,
since for this case, $p_{A}=1$. Let $P_{X^{n}}=\Unif(A)$. So, $H(X^{n})>0$,
and by the chain rule, $H(X^{n})=\sum_{i=1}^{n}H(X_{i}|X^{i-1})$.

If $P_{X^{n}Y^{n}}\in\calC_{\sym}(P_{X^{n}})$, then $Z^{n}:=X^{n}\vee Y^{n}\in A$
a.s. where $(X^{n},Y^{n})\sim P_{X^{n}Y^{n}}$. So, we have 
\[
H(Z^{n})\le\log|A|=H(X^{n}).
\]
We hence have 
\[
\sup_{P_{X^{n}Y^{n}}\in\calC_{\sym}(P_{X^{n}})}\frac{H(Z^{n})}{H(X^{n})}\le1.
\]

If $p_{A}\le t$, then $P_{X_{i}}(1)\le t,\forall i\in[n]$. Relaxing
$P_{X^{n}}=\Unif(A)$ to arbitrary distributions such that $P_{X_{i}}(1)\le t$,
we obtain $\Gamma_{n}(t)\le1$ where 
\begin{equation}
\Gamma_{n}(t):=\inf_{P_{X^{n}}:P_{X_{i}}(1)\le t,\forall i}\sup_{P_{X^{n}Y^{n}}\in\calC_{\sym}(P_{X^{n}})}\frac{H(Z^{n})}{H(X^{n})}.\label{eq:}
\end{equation}
In other words, if given $t$, $\Gamma_{n}(t)>1$, then by contradiction,
$p_{A}>t$.

We next show that $\Gamma_{n}(t)\ge\Gamma(t)$ which implies Theorem
\ref{thm:If-given-,}. To this end, we need the following lemmas.

For two conditional distributions $P_{X|U},P_{Y|V}$, denote $\mathcal{C}(P_{X|U},P_{Y|V})$
as the set of conditional distributions $Q_{XY|UV}$ such that its
marginals satisfying $Q_{X|UV}=P_{X|U},Q_{Y|UV}=P_{Y|V}$. The conditional
(Pearson) correlation coefficient of $X$ and $Y$ given $U$ is defined
by 
\[
\rho_{\mathrm{p}}(X;Y|U)=\left\{ \begin{array}{ll}
\frac{\mathbb{E}[\mathrm{cov}(X,Y|U)]}{\sqrt{\mathbb{E}[\mathrm{var}(X|U)]}\sqrt{\mathbb{E}[\mathrm{var}(Y|U)]}}, & \mathbb{E}[\mathrm{var}(X|U)]\mathbb{E}[\mathrm{var}(Y|U)]>0,\\
0, & \mathbb{E}[\mathrm{var}(X|U)]\mathbb{E}[\mathrm{var}(Y|U)]=0.
\end{array}\right.
\]
The conditional maximal correlation coefficient of $X$ and $Y$ given
$U$ is defined by 
\[
\rho_{\mathrm{m}}(X;Y|U)=\sup_{f,g}\rho_{\mathrm{p}}(f(X,U);g(Y,U)|U),
\]
where the supremum is taken over all real-valued functions $f(x,u),g(y,u)$
(such that $\mathbb{E}[\mathrm{var}(f(X,U)|U)]$, $\mathbb{E}[\mathrm{var}(g(Y,U)|U)]<\infty$).
It has been shown in \cite{yu2018conditional} that 
\[
\rho_{\mathrm{m}}(X;Y|U)=\sup_{u:P_{U}(u)>0}\rho_{\mathrm{m}}(X;Y|U=u),
\]
where $\rho_{\mathrm{m}}(X;Y|U=u)=\rho_{\mathrm{m}}(X';Y')$ with
$(X',Y')\sim P_{XY|U=u}$.
\begin{lem}[Product Construction of Couplings]
\cite[Lemma 9]{YuTan2020_exact} \cite[Corollary 3]{yu2018conditional}
\cite[Lemma 6]{beigi2015monotone} \label{lem:coupling} For any conditional
distributions $P_{X_{i}|X^{i-1}},\,P_{Y_{i}|Y^{i-1}},\,i\in[n]$ and
any 
\[
Q_{X_{i}Y_{i}|X^{i-1}Y^{i-1}}\in\mathcal{C}(P_{X_{i}|X^{i-1}},P_{Y_{i}|Y^{i-1}}),\forall i\in[n],
\]
it holds that 
\begin{align}
\prod_{i=1}^{n}Q_{X_{i}Y_{i}|X^{i-1}Y^{i-1}} & \in\mathcal{C}\Big(\prod_{i=1}^{n}P_{X_{i}|X^{i-1}},\prod_{i=1}^{n}P_{Y_{i}|Y^{i-1}}\Big).\label{eq:-15}
\end{align}
Moreover, for $(X^{n},Y^{n})\sim\prod_{i=1}^{n}Q_{X_{i}Y_{i}|X^{i-1}Y^{i-1}}$,
it holds that 
\begin{equation}
\rho_{\mathrm{m}}(X^{n};Y^{n})=\max_{i\in[n]}\rho_{\mathrm{m}}(X_{i};Y_{i}|X^{i-1},Y^{i-1}).\label{eq:-14}
\end{equation}
\end{lem}
For a conditional distribution $P_{X|U}$ defined on finite alphabets,
a conditional coupling $P_{XX'|UU'}$ of $(P_{X|U},P_{X|U})$ is called
symmetric if $P_{XX'|UU'}(x,y|u,v)=P_{XX'|UU'}(y,x|v,u)$ for all
$x,y\in\calX,u,v\in\calU$. Denote $\calC_{\sym}(P_{X|U})$ as the
set of symmetric conditional couplings of $(P_{X|U},P_{X|U})$. Applying
the lemma above to symmetric couplings, we have that if couplings
$Q_{X_{i}Y_{i}|X^{i-1}Y^{i-1}}\in\calC_{\sym}(P_{X_{i}|X^{i-1}})$
satisfy $\rho_{\mathrm{m}}(X_{i};Y_{i}|X^{i-1},Y^{i-1})\le\rho$ for
some $\rho>0$, then 
\begin{align*}
\prod_{i=1}^{n}Q_{X_{i}Y_{i}|X^{i-1}Y^{i-1}} & \in\calC_{\sym}\Big(\prod_{i=1}^{n}P_{X_{i}|X^{i-1}}\Big),\\
\rho_{\mathrm{m}}(X^{n};Y^{n}) & \le\rho,
\end{align*}
with $(X^{n},Y^{n})\sim\prod_{i=1}^{n}Q_{X_{i}Y_{i}|X^{i-1}Y^{i-1}}$.
We hence have that for any $\rho\in[0,1]$,

\begin{align}
 & \sup_{\substack{P_{X^{n}Y^{n}}\in\calC_{\sym}(P_{X^{n}}):\\
\rho_{\mathrm{m}}(X^{n};Y^{n})\le\rho
}
}H(Z^{n})\nonumber \\
 & \ge\sup_{\substack{P_{X^{n-1}Y^{n-1}}\in\calC_{\sym}(P_{X^{n-1}}):\\
\rho_{\mathrm{m}}(X^{n-1};Y^{n-1})\le\rho
}
}H(Z^{n-1})\nonumber \\
 & \qquad+\sup_{\substack{P_{X_{n}Y_{n}|X^{n-1}Y^{n-1}}\in\calC_{\sym}(P_{X_{n}|X^{n-1}}):\\
\rho_{\mathrm{m}}(X_{n};Y_{n}|X^{n-1},Y^{n-1})\le\rho
}
}H(Z_{n}|Z^{n-1})\nonumber \\
 & \ge\sup_{\substack{P_{X^{n-1}Y^{n-1}}\in\calC_{\sym}(P_{X^{n-1}}):\\
\rho_{\mathrm{m}}(X^{n-1};Y^{n-1})\le\rho
}
}H(Z^{n-1})\nonumber \\
 & \qquad+\inf_{\substack{P_{X^{n-1}Y^{n-1}}\in\calC_{\sym}(P_{X^{n-1}}):\\
\rho_{\mathrm{m}}(X^{n-1};Y^{n-1})\le\rho
}
}\sup_{\substack{P_{X_{n}Y_{n}|X^{n-1}Y^{n-1}}\in\calC_{\sym}(P_{X_{n}|X^{n-1}}):\\
\rho_{\mathrm{m}}(X_{n};Y_{n}|X^{n-1},Y^{n-1})\le\rho
}
}H(Z_{n}|Z^{n-1})\nonumber \\
 & \ge\cdots\cdots\nonumber \\
 & \ge\sum_{i=1}^{n}\inf_{\substack{P_{X^{i-1}Y^{i-1}}\in\calC_{\sym}(P_{X^{i-1}}):\\
\rho_{\mathrm{m}}(X^{i-1};Y^{i-1})\le\rho
}
}\sup_{\substack{P_{X_{i}Y_{i}|X^{i-1}Y^{i-1}}\in\calC_{\sym}(P_{X_{i}|X^{i-1}}):\\
\rho_{\mathrm{m}}(X_{i};Y_{i}|X^{i-1},Y^{i-1})\le\rho
}
}H(Z_{i}|Z^{i-1}),\label{eq:-16}
\end{align}
where the first inequality above follows by Lemma \ref{lem:coupling}
and the chain rule for entropies. {In fact, in the
derivation above, the $i$-th distribution $P_{X_{i}Y_{i}|X^{i-1}Y^{i-1}}$
is chosen as a greedy coupling in the sense that it only maximizes
the $i$-th objective function $H(Z_{i}|Z^{i-1})$, regardless of
other $H(Z_{j}|Z^{j-1})$ with $j>i$ (although it indeed affects
their values). }

By the fact that conditioning reduces entropy, it holds that 
\[
H(Z_{i}|Z^{i-1})\ge H(Z_{i}|X^{i-1},Y^{i-1}).
\]
Denote 
\begin{align}
g_{i}(P_{X^{i-1}},\rho):= & \inf_{\substack{P_{X^{i-1}Y^{i-1}}\in\calC_{\sym}(P_{X^{i-1}}):\\
\rho_{\mathrm{m}}(X^{i-1};Y^{i-1})\le\rho
}
}\sup_{\substack{P_{X_{i}Y_{i}|X^{i-1}Y^{i-1}}\in\calC_{\sym}(P_{X_{i}|X^{i-1}}):\\
\rho_{\mathrm{m}}(X_{i};Y_{i}|X^{i-1},Y^{i-1})\le\rho
}
}H(Z_{i}|X^{i-1},Y^{i-1}).\label{eq:-17-1-1-3}
\end{align}
Then, the expression at the right-hand side  of  \eqref{eq:-16} is further lower bounded by 
$\sum_{i=1}^{n}g_{i}(P_{X^{i-1}},\rho)$.  
Combing this with  \eqref{eq:} and \eqref{eq:-16}, and by
noting that $\rho\in[0,1]$ is arbitrary, we obtain that 
\begin{align}
\Gamma_{n}(t) & \ge\inf_{P_{X^{n}}:P_{X_{i}}(1)\le t,\forall i}\frac{\sup_{\rho\in[0,1]}\sum_{i=1}^{n}g_{i}(P_{X^{i-1}},\rho)}{\sum_{i=1}^{n}H(X_{i}|X^{i-1})}\nonumber \\
 & =\inf_{P_{X^{n}}:P_{X_{i}}(1)\le t,\forall i}\frac{\sup_{P_{\rho}}\mathbb{E}_{P_{\rho}}\sum_{i=1}^{n}g_{i}(P_{X^{i-1}},\rho)}{\sum_{i=1}^{n}H(X_{i}|X^{i-1})}\nonumber \\
 & \ge\sup_{P_{\rho}}\inf_{P_{X^{n}}:P_{X_{i}}(1)\le t,\forall i}\frac{\sum_{i=1}^{n}\mathbb{E}_{P_{\rho}}g_{i}(P_{X^{i-1}},\rho)}{\sum_{i=1}^{n}H(X_{i}|X^{i-1})}\nonumber \\
 & \ge\sup_{P_{\rho}}\inf_{P_{X^{n}}:P_{X_{i}}(1)\le t,\forall i}\min_{i\in[n]:H(X_{i}|X^{i-1})>0}\frac{\mathbb{E}_{P_{\rho}}g_{i}(P_{X^{i-1}},\rho)}{H(X_{i}|X^{i-1})}\label{eq:-3}\\
 & \ge\sup_{P_{\rho}}\inf_{P_{X^{j}}:H(X_{j}|X^{j-1})>0,P_{X_{j}}(1)\le t}\frac{\mathbb{E}_{P_{\rho}}g_{j}(P_{X^{j-1}},\rho)}{H(X_{j}|X^{j-1})},\nonumber 
\end{align}
where
\begin{itemize}
\item \eqref{eq:-3} follows since $\frac{a+b}{c+d}\ge\min\{\frac{a}{c},\frac{b}{d}\}$
for $a,b\ge0,c,d>0$, and $H(X_{i}|X^{i-1})=0$ implies $X_{i}$ is
a function of $X^{i-1}$, and hence, $g_{i}(P_{X^{i-1}},\rho)=0$;
\item the index $j$ in the last line is the optimal $i$ attaining the minimum
in \eqref{eq:-3}.
\end{itemize}
Denote $X=X_{j},Y=Y_{j},U=X^{j-1},V=Y^{j-1}$, and $Z=X\lor Y$. Then,
\begin{align}
\Gamma_{n}(t) & \ge\sup_{P_{\rho}}\inf_{P_{UX}:H(X|U)>0,P_{X}(1)\le t}\mathbb{E}_{P_{\rho}}\left[\inf_{\substack{P_{UV}\in\calC_{\sym}(P_{U}):\\
\rho_{\mathrm{m}}(U;V)\le\rho
}
}\sup_{\substack{P_{XY|UV}\in\calC_{\sym}(P_{X|U}):\\
\rho_{\mathrm{m}}(X;Y|U,V)\le\rho
}
}\frac{H(Z|U,V)}{H(X|U)}\right].\label{eq:-2}
\end{align}

We next further simplify the lower bound in \eqref{eq:-2}. Denote
\begin{equation}
p=P_{X|U}(1|U),q=P_{Y|V}(1|V),r=P_{XY|UV}(1,1|U,V).\label{eq:pq}
\end{equation}
So, 
\[
P_{XY|UV}(\cdot|U,V)=\begin{bmatrix}1+r-p-q & q-r\\
p-r & r
\end{bmatrix}
\]
with 
\[
\max\{0,p+q-1\}\le r\le\min\{p,q\}.
\]
Note that 
\begin{align}
\rho_{\mathrm{m}}(X;Y|U,V) & =\sup_{u,v:P_{UV}(u,v)>0}\rho_{\mathrm{m}}(X_{u};Y_{v})\nonumber \\
 & =\sup_{u,v:P_{UV}(u,v)>0}\left|\rho_{\mathrm{p}}(X_{u};Y_{v})\right|\label{eq:-1}\\
 & =\sup_{u,v:P_{UV}(u,v)>0}\frac{\left|r-pq\right|}{\sqrt{p(1-p)q(1-q)}},\nonumber 
\end{align}
where $\left(X_{u},Y_{v}\right)\sim P_{XY|U=u,V=v}$, $\rho_{\mathrm{p}}$
denotes the Pearson correlation coefficient, and \eqref{eq:-1} follows
since the maximal correlation coefficient between two binary random
variables is equal to the absolute value of the Pearson correlation
coefficient between them; see, e.g., \cite{anantharam2013maximal}.
So, $\rho_{\mathrm{m}}(X;Y|U,V)\le\rho$ is equivalent to $\frac{\left|r-pq\right|}{\sqrt{p(1-p)q(1-q)}}\le\rho$
a.s., and also equivalent to $z_{1}\le r\le z_{2}$ a.s.

The inner supremum in \eqref{eq:-2} can be rewritten as 
\begin{align*}
 & \sup_{P_{XY|UV}\in\calC_{\sym}(P_{X|U}):\rho_{\mathrm{m}}(X;Y|U,V)\le\rho}H(Z|U,V)\\
 & =\mathbb{E}_{p,q}\sup_{\max\{0,p+q-1,z_{1}\}\le r\le\min\{p,q,z_{2}\}}h(p+q-r).
\end{align*}
By the fact that $h$ is increasing on $[0,1/2]$ and decreasing on
$[1/2,1]$, it holds that the optimal $r$ attaining the supremum
in the last line above, denoted by $r^{*}$, is the median of $\max\{0,p+q-1,z_{1}\}$,
$p+q-1/2$, and $\min\{p,q,z_{2}\}$, which implies 
\[
p+q-r^{*}=\varphi(\rho,p,q).
\]
Recall the definition of $\varphi$ in \eqref{eq:-4}. So, the inner
supremum in \eqref{eq:-2} is equal to $\frac{\mathbb{E}_{p,q}h(\varphi(\rho,p,q))}{\mathbb{E}h(p)}$. 

{We make following observations. Firstly, 
\begin{align*}
H(X|U) & =\mathbb{E}h(p),\\
P_{X}(1) & =\mathbb{E}p.
\end{align*}
Secondly, by the definition of maximal correlation, $\rho_{\mathrm{m}}(p;q)\le\rho_{\mathrm{m}}(U;V)$
holds (which is known as the data processing inequality) since $p,q$
are respectively functions of $U,V$; see \eqref{eq:pq}. Lastly,
observe that $P_{UV}$ is symmetric, and $p,q$ are obtained from
$U,V$ via the same function $P_{X|U}(1|\cdot)$ (since $P_{X|U}=P_{Y|V}$
holds by the symmetry of $P_{XY|UV}$). Hence, $P_{pq}$ is symmetric
as well.  Substituting all of these into \eqref{eq:-2} yields $\Gamma_{n}(t)\ge\Gamma(t)$.}

\section{Proof of Proposition \ref{prop:For-,where-the}}

By choosing $P_{\rho}=(1-\alpha)\delta_{0}+\alpha\delta_{1}$ in \eqref{eq:-5},
we obtain 
\[
\Gamma(t)\ge\sup_{\alpha\in[0,1]}\inf_{\textrm{symmetric }P_{pq}:\mathbb{E}h(p)>0,\mathbb{E}p\le t}\frac{g(P_{pq},\alpha)}{\mathbb{E}h(p)}.
\]
Note that $P_{pq}\mapsto g(P_{pq},\alpha)$ is concave, since by \cite[Lemma 5]{alweiss2022improved},
$P_{p}\mapsto\mathbb{E}_{(p,q)\sim P_{p}^{\otimes2}}h(p+q-pq)$ is
concave, and $P_{pq}\mapsto P_{p}$ is linear.

Let $B$ be a finite subset of $[0,1]$. Let $\mathcal{P}_{B}$ be
the set of symmetric distributions $P_{pq}$ concentrated on $B^{2}$
such that $\mathbb{E}p\le t$. By the Krein--Milman theorem, $\mathcal{P}_{B}$
is equal to the closed convex hull of its extreme points. These extreme
points are of the form $(1-\beta)Q_{a_{1},a_{2}}+\beta Q_{b_{1},b_{2}}$
with $0\le a\le t<b\le1$ and $\beta=0$ or $\frac{t-a}{b-a}$, where
recall the definitions $a:=\frac{a_{1}+a_{2}}{2},b:=\frac{b_{1}+b_{2}}{2}$,
and $Q_{x,y}:=\frac{1}{2}\delta_{(x,y)}+\frac{1}{2}\delta_{(y,x)}$
in \eqref{eq:ab} and \eqref{eq:Q}. By Carath\'eodory's theorem, it
is easy to see that the convex hull of these extreme points is closed
(in the weak topology, or equivalently, in the relative topology on
the probability simplex). So, every $P_{pq}$ supported on a finite
set $B^{2}\subseteq[0,1]^{2}$ such that $\mathbb{E}p\le t$ is a
convex combination of the extreme points above, i.e., $P_{pq}=\sum_{i=1}^{k}\gamma_{i}Q_{i}$
where $Q_{i},i\in[k]$ are extreme points, and $\gamma_{i}>0$ and
$\sum_{i=1}^{k}\gamma_{i}=1$. For this distribution, 
\begin{align*}
\frac{g(P_{pq},\alpha)}{\mathbb{E}h(p)} & =\frac{g(\sum_{i=1}^{k}\gamma_{i}Q_{i},\alpha)}{\sum_{i=1}^{k}\gamma_{i}\mathbb{E}_{Q_{i}}h(p)}\\
 & \ge\frac{\sum_{i=1}^{k}\gamma_{i}g(Q_{i},\alpha)}{\sum_{i=1}^{k}\gamma_{i}\mathbb{E}_{Q_{i}}h(p)}\\
 & \ge\min_{i:\mathbb{E}_{Q_{i}}h(p)>0}\frac{g(Q_{i},\alpha)}{\mathbb{E}_{Q_{i}}h(p)}
\end{align*}
where in the last line, we use the fact that $\mathbb{E}_{Q_{i}}h(p)=0$
implies $Q_{i}=\delta_{(0,0)}$ (note that $t<1/2$), and hence, $g(Q_{i},\alpha)=0$.

Therefore, 
\begin{equation}
\Gamma(t)\ge\sup_{\alpha\in[0,1]}\inf_{P_{pq}:\mathbb{E}h(p)>0}\frac{g(P_{pq},\alpha)}{\mathbb{E}h(p)},\label{eq:-8}
\end{equation}
where the infimum is taken over distributions $P_{pq}$ of the form
$(1-\beta)Q_{a_{1},a_{2}}+\beta Q_{b_{1},b_{2}}$ with $0\le a\le t<b\le1$
and $\beta=0$ or $\beta=\frac{t-a}{b-a}>0$ such that $\mathbb{E}h(p)>0$.
(Recall the definition of $a,b$ in \eqref{eq:ab}.)

\section{Discussion}

The breakthrough made by Gilmer \cite{gilmer2022constant} shows the
power of information-theoretic techniques in tackling problems from
related fields. In fact, the union-closed sets conjecture has a natural
interpretation in the information-theoretic (or coding-theoretic)
sense. Consider the memoryless OR multi-access channel $(x^{n},y^{n})\in\Omega^{2n}\mapsto x^{n}\lor y^{n}\in\Omega^{n}$.
We would like to find a nonempty code $A\subseteq\Omega^{n}$ to generate
two independent inputs $X^{n},Y^{n}$ with each following $\Unif(A)$
such that the input constraint $\mathbb{E}[X_{i}]\le t,\forall i\in[n]$
is satisfied and the output $X^{n}\lor Y^{n}$ is still in $A$ a.s.
The union-closed sets conjecture states that such a code exists if
and only if $t\ge1/2$. Based on this information-theoretic interpretation,
it is reasonable to see that the information-theoretic techniques
work for this conjecture. It is well-known that information-theoretic
techniques usually work very well for problems with ``approximate''
constraints, e.g., the channel coding problem with the asymptotically
vanishing error probability constraint (or the approximate version
of the union-closed sets problem introduced in \cite{chase2022approximate}).
It is hard to say whether information-theoretic techniques are sufficient
to prove sharp bounds for problems with ``exact'' constraints, e.g.,
the zero-error coding problem (or the original version of the union-closed
sets conjecture).

Furthermore, as an intermediate result, it has been shown that $\Gamma_{n}(t)>1$
implies $p_{A}>t$ for any OR-closed $A\subseteq\Omega^{n}$. Here
$\Gamma_{n}(t)$ is given in \eqref{eq:}, expressed in the multi-letter
form {(i.e., the dimension-dependent form)}. By the
super-block coding argument, it is verified that given $t>0$,
$\lim_{n\to\infty}\Gamma_{n}(t)$ exists. It is interesting to investigate
this limit, and prove a single-letter {(dimension-independent)}
expression for it.

For simplicity, in this paper, we only consider the maximal correlation
coefficient as the constraint function. In fact, the maximal correlation
coefficient used here can be replaced by other functionals. The key
property of the maximal correlation coefficient we used in this paper
is the ``tensorization'' property, i.e., \eqref{eq:-14} (in fact,
only ``$\le$'' part of \eqref{eq:-14} was used in our proof).
In literature, there is a class of measures of correlation satisfying
this property, e.g., the hypercontractivity constant, strong data
processing inequality constant, or more generally, $\Phi$-ribbons,
see \cite{ahlswede1976spreading,raginsky2016strong,beigi2018phi}.
(Although the tensorization property in the literature is only defined
and proven for independent   random variables, this property can
be   extended to the coupling constructed in \eqref{eq:-15}.)
Following the same proof steps given in this paper, one can  
obtain various variants of Theorem \ref{thm:If-given-,} with the
maximal correlation coefficient replaced by other quantities, as long
as these quantities satisfy the tensorization property. Another potential
direction is to replace the Shannon entropy with a class of more general
quantities, R\'enyi entropies. However, unfortunately   R\'enyi entropies   do  not satisfy the chain rule (unlike the Shannon entropy), which leads to   a serious difficulty   in  single-letterizing the corresponding  multi-letter bound  like $\Gamma_{n}(t)$ in \eqref{eq:} {(i.e.,
in making the multi-letter bound dimension-independent)}.

\vspace{6pt}




\funding{This work was supported by the NSFC grant 62101286 and
the Fundamental Research Funds for the Central Universities of China
(Nankai University).}

\institutionalreview{Not applicable.}

\informedconsent{Not applicable.}

\dataavailability{Not applicable.}

\acknowledgments{The author would like to thank Fan Chang for bringing
Gilmer's breakthrough \cite{gilmer2022constant} to his attention, and thank Stijn Cambie for sharing his early draft of \cite{cambie2022better}.
The author also would like to thank the guest editor, Prof. Igal Sason,
for his invitation to submit this paper to the Entropy, and thank him and 
the anonymous referees for their comments, which led to significant
improvements in the presentation of this paper.}

\conflictsofinterest{Not applicable.}

\begin{adjustwidth}{-\extralength}{0cm} 

\reftitle{References}




\bibliographystyle{plain}
\bibliography{ref}

\PublishersNote{} \end{adjustwidth}
\end{document}